\documentclass{amsart}
\begin{document}
\title[Differential cocycles]{Arithmetic differential equations on $GL_n$, I:\\
differential cocycles}
\author{Alexandru Buium and Taylor Dupuy}
\def \Rp{R_p}
\def \Rpi{R_{\pi}}
\def \dpi{\d_{\pi}}
\def \bT{{\bf T}}
\def \cI{{\mathcal I}}
\def \cH{{\mathcal H}}
\def \cJ{{\mathcal J}}
\def \ZN{\bZ[1/N,\zeta_N]}
\def \tA{\tilde{A}}
\def \o{\omega}
\def \tB{\tilde{B}}
\def \tC{\tilde{C}}
\def \alph{A}
\def \bet{B}
\def \bsigma{\bar{\sigma}}
\def \y{^{\infty}}
\def \Ra{\Rightarrow}
\def \uBS{\overline{BS}}
\def \lBS{\underline{BS}}
\def \lB{\underline{B}}
\def \<{\langle}
\def \>{\rangle}
\def \hL{\hat{L}}
\def \cU{\mathcal U}
\def \cF{\mathcal F}
\def \S{\Sigma}
\def \st{\stackrel}
\def \sd{Spec_{\d}\ }
\def \pd{Proj_{\d}\ }
\def \s{\sigma_2}
\def \i{\sigma_1}
\def \bs{\bigskip}
\def \cD{\mathcal D}
\def \cC{\mathcal C}
\def \cT{\mathcal T}
\def \cK{\mathcal K}
\def \cX{\mathcal X}
\def \sX{X_{set}}
\def \cY{\mathcal Y}
\def \cS{X}
\def \cR{\mathcal R}
\def \cE{\mathcal E}
\def \tcE{\tilde{\mathcal E}}
\def \cP{\mathcal P}
\def \cA{\mathcal A}
\def \cV{\mathcal V}
\def \cM{\mathcal M}
\def \cL{\mathcal L}
\def \cN{\mathcal N}
\def \tcM{\tilde{\mathcal M}}
\def \caS{\mathcal S}
\def \cG{\mathcal G}
\def \cB{\mathcal B}
\def \tG{\tilde{G}}
\def \cF{\mathcal F}
\def \h{\hat{\ }}
\def \hp{\hat{\ }}
\def \tS{\tilde{S}}
\def \tP{\tilde{P}}
\def \tA{\tilde{A}}
\def \tX{\tilde{X}}
\def \tcS{\tilde{X}}
\def \tT{\tilde{T}}
\def \tE{\tilde{E}}
\def \tV{\tilde{V}}
\def \tC{\tilde{C}}
\def \tI{\tilde{I}}
\def \tU{\tilde{U}}
\def \tG{\tilde{G}}
\def \tu{\tilde{u}}
\def \chu{\check{u}}
\def \tx{\tilde{x}}
\def \tL{\tilde{L}}
\def \tY{\tilde{Y}}
\def \d{\delta}
\def \e{\chi}
\def \bW{\mathbb W}
\def \bV{{\mathbb V}}
\def \bF{{\bf F}}
\def \bE{{\bf E}}
\def \bC{{\bf C}}
\def \bO{{\bf O}}
\def \bR{{\bf R}}
\def \bA{{\bf A}}
\def \bB{{\bf B}}
\def \cO{\mathcal O}
\def \ra{\rightarrow}
\def \bx{{\bf x}}
\def \f{{\bf f}}
\def \bX{{\bf X}}
\def \bH{{\bf H}}
\def \bS{{\bf S}}
\def \bF{{\bf F}}
\def \bN{{\bf N}}
\def \bK{{\bf K}}
\def \bE{{\bf E}}
\def \bB{{\bf B}}
\def \bQ{{\bf Q}}
\def \bd{{\bf d}}
\def \bY{{\bf Y}}
\def \bU{{\bf U}}
\def \bL{{\bf L}}
\def \bQ{{\bf Q}}
\def \bP{{\bf P}}
\def \bR{{\bf R}}
\def \bC{{\bf C}}
\def \bD{{\bf D}}
\def \bM{{\bf M}}
\def \bZ{{\mathbb Z}}
\def \xtoleqr{x^{(\leq r)}}
\def \hU{\hat{U}}
\def \k{\kappa}
\def \ee{\overline{p^{\k}}}

\newtheorem{THM}{{\!}}[section]
\newtheorem{THMX}{{\!}}
\renewcommand{\theTHMX}{}
\newtheorem{theorem}{Theorem}[section]
\newtheorem{corollary}[theorem]{Corollary}
\newtheorem{lemma}[theorem]{Lemma}
\newtheorem{proposition}[theorem]{Proposition}
\theoremstyle{definition}
\newtheorem{definition}[theorem]{Definition}
\theoremstyle{remark}
\newtheorem{remark}[theorem]{Remark}
\newtheorem{example}[theorem]{\bf Example}
\numberwithin{equation}{section}
\address{University of New Mexico \\ Albuquerque, NM 87131}
\email{buium@math.unm.edu} 
%\subjclass[2000]{....}
\maketitle

\begin{abstract}
The theory of differential equations has an arithmetic analogue \cite{book} in which 
derivatives are replaced by Fermat quotients. One can then ask what is the arithmetic analogue of a linear differential equation. The study of usual linear differential equations is the same as the study 
of the differential cocycle   from $GL_n$ into its Lie algebra given by the logarithmic derivative \cite{kolchin}. However we prove here that there are no such cocycles in the context of arithmetic differential equations. In sequels of this paper \cite{adel2, adel3} we will remedy the situation by introducing arithmetic analogues of Lie algebras and a skew version of differential cocycles; this will lead to a theory of linear arithmetic differential equations.
\end{abstract}

\section{Introduction and main results}

In \cite{char} an arithmetic analogue of differential equations was introduced in which derivations are replaced by Fermat quotient operators; cf. \cite{book} for an overview of the theory which was mostly concerned with Abelian varieties and Shimura varieties. With the exception of \cite{simple}, however, little attention has been given to the case of arithmetic differential equations attached to linear algebraic groups such as $GL_n$. The present paper is the first in a series of papers whose  purpose is to shed some light into the theory for $GL_n$.
In particular, one of our main motivations is to understand  what is the correct notion of
 ``linear arithmetic differential equation".
 
 The basic clue should come from the Ritt-Kolchin differential algebra \cite{kolchin} whose arithmetic analogue is the theory in \cite{char,book}. So let us examine Kolchin's setting first.
 Since we will later treat Kolchin's theory and our arithmetic theory simultaneously it is convenient to use the same notation in two different contexts. We call these contexts the {\it $\delta$-algebraic} setting (corresponding to Kolchin's theory) and the {\it  $\delta$-arithmetic} setting (corresponding to the theory in \cite{char, book}). 
 
 In the $\delta$-algebraic setting we denote by $R$
 a field of characteristic zero equipped with a derivation $\d:R\ra R$ and we assume $R$ is $\d$-{\it closed} (which is the same as {\it constrainedly closed} in Kolchin's terminology \cite{kolchin}). We denote by $R^{\d}$ the field of constants $\{c\in R;\d c=0\}$. In this context we will consider smooth schemes of finite type $X$ over $R$ (i.e. nonsingular varieties) and we denote by $X(R)$ the set of $R$-points of $X$; if there is no danger of confusion we often simply write $X$ in place of $X(R)$. 
 
 In the  $\delta$-arithmetic setting we assume $R$ is the unique complete discrete valuation ring with maximal ideal generated by an odd prime $p$ and  residue field $k=R/pR$ equal to the algebraic closure ${\mathbb F}_p^a$ of ${\mathbb F}_p$; then we denote by $\d:R\ra R$ the unique $p$-derivation on $R$ in the sense of \cite{char}; recall that $\d x=\frac{\phi(x)-x^p}{p}$ where $\phi:R\ra R$ is the unique ring homomorphism lifting the $p$-power Frobenius on the residue field $k$. We denote by 
$R^{\d}$ the monoid of constants $\{\lambda\in R;\d \lambda=0\}$; so $R^{\d}$ consists of $0$ and all roots of unity in $R$. Also we denote by $K$ the fraction field of $R$.
 In this context we will consider smooth schemes of finite type $X$ over $R$ or, more generally,  smooth $p$-formal schemes of finite type, by which we mean formal schemes locally isomorphic to $p$-adic completions of smooth schemes of finite type; we denote by $X(R)$ the set of $R$-points of $X$; if there is no danger of confusion we often simply write $X$ in place of $X(R)$. Groups in the category of smooth $p$-formal schemes will be called smooth group $p$-formal schemes.
 
In the $\d$-algebraic setting a map $f:R^N\ra R^M$
will be called a $\d$-map of order $n$ if there exists an $M$-vector  $F=(F_j)$ of polynomials $F_j\in R[x_0,...,x_n]$ with coefficients in $R$ in $n+1$ $N$-tuples of variables $x_0,...,x_n$, such that 
\begin{equation}
\label{singe}
f(a)=F(a,\d a,...,\d^n a),\ \ a\in R^N.\end{equation}

In the $\d$-arithmetic setting a map $f:R^N\ra R^M$
will be called a $\d$-map of order $n$ if there exists an $M$-vector  $F=(F_j)$ of restricted power series $F_j\in R[x_0,...,x_n]\h$ (where $\h$ means $p$-adic completion) such that \ref{singe} holds.
 In both the $\d$-algebraic and the $\d$-arithmetic setting one can then consider affine smooth schemes $X,Y$ and define a map $f:X(R)\ra Y(R)$ to be a $\d$-map of order $n$ if there exist  embeddings $X\subset {\mathbb A}^N$, $Y\subset {\mathbb A}^M$ such that $f$ is induced by a $\d$-map $R^N\ra R^M$ of order $n$; we simply write $f:X\ra Y$. Even more generally if $X$ is any smooth scheme and $Y$ is an affine scheme a set theoretic map $X\ra Y$ is called a $\d$-map of order $n$ if there exists an affine cover $X=\bigcup X_i$ such that all the induced maps $X_i\ra Y$ are $\d$-maps of order $n$.
 If $G,H$ are smooth group schemes a $\d$-homomorphism $G\ra H$ is, by definition,  a $\d$-map which is also a homomorphism. 
 We shall review these concepts in section 2 following 
\cite{char,book} in the  $\delta$-arithmetic setting and 
 \cite{kolchin, cassidy, hermann} in the $\delta$-algebraic setting. We will use a slightly different (but equivalent) approach using jet spaces.
  
    Assume we are in either the $\delta$-algebraic or in the $\delta$-arithmetic  setting.  Let $G$ be an algebraic group (i.e. smooth group scheme) over $R$ and let $L(G)$ be its Lie algebra which we view as a group scheme over $R$ isomorphic as a scheme with the affine space. Denote by
  $\star:G\times L(G)\ra L(G)$ the (left) adjoint action.
  
  \begin{definition}
   A classical $\d$-cocycle of $G$ with values in $L(G)$ is  a $\d$-map $f:G\ra L(G)$ which is a cocycle for the adjoint action; i.e. for all $g_1,g_2\in G(R)$ we have
 \begin{equation}
 \label{theone}
 f(g_1g_2)=f(g_1)+g_1\star f(g_2).\end{equation}\end{definition}

 Here $G(R)$ is the set of $R$-points of $G$.
 If in the definition above $f$ is a regular map (i.e. a morphism of schemes over $R$) then we say $f$ is a regular cocycle. 
 
 The adjective {\it classical} in the above definition was used in order to distinguish between the {\it classical $\d$-cocycles} introduced above and {\it skew $\d$-cocycles} to be introduced in \cite{adel2}; skew $\d$-cocycles will then   be viewed as ``non-classical" objects.
 
 \begin{remark}\ 
 
 1) For $G$ a smooth closed subgroup scheme of some $GL_n$, identifying $L(G)$ with a Lie subalgebra of ${\mathfrak gl}_n=L(GL_n)$, we get that  condition \ref{theone} reads
 $$f(g_1g_2)=f(g_1)+g_1f(g_2)g_1^{-1}.$$
 On the other hand if $G$ is a (non-necessarily affine) smooth commutative group then  condition
 \ref{theone} reads
 $$f(g_1g_2)=f(g_1)+f(g_2)$$
 i.e. $f$ is simply a $\d$-homomorphism.
 In view of Kolchin's theory \cite{kolchin} it seems reasonable to introduce the following definition  (in both the $\delta$-algebraic and the  $\delta$-arithmetic contexts): a linear differential equation
 is an equation of the form
 $f(u)=\alpha$
 where $f$ is a classical $\d$-cocycle of $G$ with values in $L(G)$, $u\in G(R)$ is the unknown and $\alpha \in L(G)(R)$ is a given $R$-point of the Lie algebra. 
 
 2)  Assume we are in the $\delta$-algebraic context and assume $G$ is defined over the field of constants $R^{\d}$ of $\d$. In that case  there is a remarkable order $1$ classical $\d$-cocycle $l\d:G\ra L(G)$ called the  {\it Kolchin logarithmic derivative}.
 For $G$ closed in $GL_n$ and $\d$-horizontal (i.e. defined by equations with coefficients in $R^{\d}$) $l\d$ is defined by $l\d g=\d g \cdot g^{-1}$, $g\in G(R)$. This essentially models, in the context at hand, concepts that go back to Lie and Cartan. Also the equation $l\d(u)=\alpha$ for $G=GL_n$ reduces to an equation $\d u=\alpha u$ which is the familiar form of a  linear differential equation. The above construction can be generalized to the case
  when $G$ has a structure of $D$-group \cite{lnm,pillay} i.e. $G$ comes equipped with a derivation $\cO_G\ra \cO_G$ extending that of $R$ and compatible with the group structure.
 
 3) Assume again that we are in the $\d$-algebraic context. Then
 for any Abelian variety $A$ over $R$ (not necessarily defined over the field of constants $R^{\d}$!) there is a surjective, order $2$,  $\d$-homomorphisms $f:A\ra L(A)$
 (hence a classical $\d$-cocycle) into the Lie algebra $L(A)$ of $A$; this follows from Manin's work \cite{manin} (cf. \cite{annals,hermann} for a different construction of Manin's map $f$). It turns out that there is an analogue of Manin maps in the  $\delta$-arithmetic setting \cite{char}: for any abelian scheme $A$ over $R$ there is an order $2$ surjective $\d$-homomorphism $f:A\ra L(A)$ with remarkable properties (e.g one has a precise description of its kernel which is analogous to the Manin theorem of the kernel \cite{manin}). 
 Similar maps exist for linear tori in place of abelian schemes; in particular there is a complete description  of the $\d$-homomorphisms ${\mathbb G}_m\ra {\mathbb G}_a=L({\mathbb G}_m)$;
cf. Lemma \ref{allchar} (which follows immediately from \cite{char}).
 \end{remark}
 
 To summarize, one has a good definition of linear differential equations in the $\delta$-algebraic setting. Moreover there is a nice $\d$-arithmetic analogue of linear differential equations in the case of Abelian varieties.  However our Theorem
 \ref{main}  below shows that there is no ``naive" $\d$-arithmetic analogue of the 
 linear differential equations in the case of $GL_n$: what we show is that there is no ``naive" $\d$-arithmetic analogue of the 
 Kolchin logarithmic derivative for $GL_n$, $n\geq 2$.

\begin{theorem}\label{main}
Assume we are in the  $\delta$-arithmetic  setting.
Let $f:GL_n\ra {\mathfrak gl}_n$ be a classical $\d$-cocycle. Then there exists a $\d$-ho\-mo\-morphism $\omega:{\mathbb G}_m\ra {\mathbb G}_a$ and there exists $v\in gl_n(R)$ such that
 for all $g\in GL_n(R)$ we have: 
$$f(g)=\omega(\det(g)) 1_n +gvg^{-1}-v.$$\end{theorem}

The  corresponding statement in algebraic geometry, saying that if a  map of varieties
$GL_n\ra {\mathfrak gl}_n$ over an algebraically closed field of characteristic zero
is  a cocycle then the map must be a couboundary, 
is ``well known" and  follows easily from Whitehead's lemma; cf. Lemma \ref{cocoa}.

Note also that our computations have a $\delta$-algebraic variant leading to a characterization of Kolchin's logarithmic  derivative that seems to be new. Indeed assume we are in the $\delta$-algebraic setting and let $f:GL_n\ra {\mathfrak gl}_n$ be a classical $\d$-cocycle. Let us say that $f$ is $\d$-coherent if for any algebraic subgroup $G\subset GL_n$ defined over the field of constants $R^{\d}$ with Lie algebra $L(G)$ we have that $f(G)\subset L(G)$. It is well known (and trivial to check) that Kolchin's logarithmic derivative $l\d:GL_n\ra {\mathfrak gl}_n$, $l\d g=\d g \cdot g^{-1}$,  is a $\d$-coherent classical $\d$-cocycle.
Conversely we will prove:

\begin{theorem}
\label{secondary}
Assume we are in the $\delta$-algebraic setting.
Let $f:GL_n\ra {\mathfrak gl}_n$ be a $\d$-coherent classical $\d$-cocycle and assume $n\geq 2$. Then there exists $\nu\in R$ such that for all $g\in GL_n(R)$,
$$f(g)=\nu \cdot \d g \cdot g^{-1}.$$
\end{theorem}

For $n=1$ the above fails; instead one has a complete description of the situation in this case due to Cassidy \cite{cassidy}; cf. Lemma \ref{allchar} for a review of this.

\medskip

As already mentioned Theorem \ref{main} shows that there is no ``naive" analogue of Kolchin's logarithmic derivative for $GL_n$  when $n\geq 2$. In order to find, then, an arithmetic analogue of linear differential equations we will be led in \cite{adel2, adel3} to replace Lie algebras of algebraic groups by objects that are less linear and, arguably,  are better adapted to the arithmetic jet theory. 
This will lead to a theory of linear arithmetic differential equations. 

The present paper is organized as follows. Section 2 reviews some of the basic concepts 
in \cite{char,book}. In section 3 we present the proofs of Theorems \ref{main} and \ref{secondary} about classical $\d$-cocycles. 

\medskip

{\bf Acknowledgement}.
The authors are indebted to P. Cartier for inspiring discussions and to A. Minchenko for providing the proof of Lemma \ref{cocoa}. Also the first author would like to acknowledge partial support from the Hausdorff Institute of Mathematics in Bonn
and from the NSF through grant DMS 0852591.

\section{Review of $p$-jets}

The aim of this section is to review the relevant material on the $\delta$-arithmetic setting in \cite{char,book}. We will also include the relevant corresponding comments on differential algebra \cite{kolchin, cassidy, hermann}. 

Assume first we are in the $\delta$-arithmetic setting. 
Recall from the introduction that we denote by $R$  the complete discrete valuation ring with maximal ideal generated by an odd prime $p$ and algebraically closed residue field $k={\mathbb F}_p^a$. Then $R$ comes equipped with a unique lift of Frobenius $\phi:R\ra R$ i.e. with a ring homomorphism whose reduction mod $p$ is the $p$-th power map. 
For 
$x$ a tuple of indeterminates over $R$ and tuples of indeterminates $x',...,x^{(n)},...$
we let $R\{x\}=R[x,x',x'',...]$ and we still denote by $\phi:R\{x\}\ra R\{x\}$  the unique lift of Frobenius
extending $\phi$ on $R$ such that $\phi(x)=x^p+px'$, $\phi(x')=(x')^p+px''$, etc. Then we let
 $\d:R\{x\}\ra R\{x\}$ be defined as $\d f=p^{-1}(\phi(f)-f^p)$; so $\d x=x'$, $\d x'=x''$, etc.  We view $\d$ as an analogue of the total derivative operator in differential algebra. 
 Now
for any affine scheme of finite type $X=Spec\ R[x]/(f)$ over $R$, where $f$ is a tuple of polynomials, we define the $p$-jet spaces of $X$ as being the $p$-formal schemes
\begin{equation}
\label{fort}
J^n(X)=Spf\ \frac{R[x,x',...,x^{(n)}]\h}{(f,\d f,...,\d^n f)}\end{equation}
For $X$ of finite type but not necessarily affine we define $J^n(X)=\bigcup J^n(X_i)$ where $X=\bigcup X_i$ is an affine cover and the gluing is an obvious one. The spaces $J^n(X)$ have an obvious universality property for which we refer to \cite{char,book} and can be defined for $X$ a $p$-formal scheme  of finite type as well. If $X/R$ is smooth then $J^n(X)$ is locally the $p$-adic completion of a smooth scheme.
The universality property yields natural maps on sets of $R$-points $\nabla^n:X(R)\ra J^n(X)(R)$; for $X={\mathbb A}^m$ (the affine space),
$J^n(X)={\mathbb A}^{m(n+1)}$ and $\nabla^n(a)=(a,\d a,...,\d^na)$.
Let $X,Y$ be schemes of finite type over $R$; by a $\d$-map (of order $n$) $f:X\ra Y$ we understand a map of $p$-formal schemes $J^n(X)\ra J^0(Y)=\widehat{Y}$. Two $\d$-maps 
$X\ra Y$ and $Y\ra Z$ of orders $n$ and $m$ respectively can be composed (using the universality property) to yield a $\d$-map of order $n+m$. Any $\d$-map $f:X\ra Y$ induces a set theoretic map  $f_*:X(R)\ra Y(R)$ defined by $f_*(P)=f(\nabla^n(P))$; if $X,Y$ are smooth the map $f_*$ uniquely determines the map $f$ and, in this case, we simply write $f$ instead of $f_*$ (and $X,Y$ instead of $X(R),Y(R)$).
A $\d$-map $X\ra Y$ of order zero 
is nothing but a map of $p$-formal schemes $\widehat{X}\ra \widehat{Y}$.
The functors $J^n$ commute with products and send groups into groups. By a $\d$-homomorphism
$f:G\ra H$ between two group schemes (or group $p$-formal schemes) we understand a group homomorphism $J^n(G)\ra J^0(H)=\widehat{H}$. 

The Kolchin $\delta$-algebraic setting can be presented in a similar way. In this setting one defines jet spaces $J^n(X)$ of schemes of finite type $X$ over $R$ by the same formula \ref{fort} in which one drops the symbol $\h$, one replaces $Spf$ by $Spec$, and one takes $\d:R\{x\}\ra R\{x\}$ to be the unique derivation extending $\d$ on $R$ and sending $\d x=x'$, $\d x'=x''$, etc. If $X$ is smooth then $J^n(X)$ are  smooth over $R$. Again one has natural set theoretic maps $\nabla^n:X(R)\ra J^n(X)(R)$.  A $\d$-map $f:X\ra Y$ is, as before,  a morphism of schemes $J^n(X)\ra Y$; one defines similarly the composition of $\d$-maps. Any $\d$-map $f:X\ra Y$ induces a set theoretic map $f_*:X(R)\ra Y(R)$ which determines $f$ uniquely if $X$ and $Y$ are smooth; in this case we write $f$ in place of $f_*$. A $\d$-map of order zero is nothing but a map of schemes. One defines $\d$-homomorphisms in the expected way. The formalism of flows and prime integrals is similar and boils down, of course, to the classical picture.

\begin{remark}
The $p$-jet spaces $J^n(X)$ have natural algebraizations which, in their turn, have multi-prime versions. The resulting jet functors are adjoint to the Witt functors. Cf. Borger's work \cite{borger, borger2}. 
\end{remark}

\section{Proofs of Theorems \ref{main} and \ref{secondary}}

We start with the ``non-differential" case which is
   ``well known"; for convenience we include a proof (for which we are indebted to A. Minchenko):

\begin{lemma}\label{cocoa}
Let $f:GL_n(F) \ra {\mathfrak gl}_n(F)$ be a regular map of algebraic varieties over an algebraically closed field $F$ of characteristic zero and assume $f$ is a cocycle  i.e. for all $g_1,g_2\in GL_n(F)$ equation \ref{theone} holds. Then 
 $f$ 
is a coboundary, i.e. there exists $v\in gl_n(F)$ such that for all $g\in GL_n(F)$,
 $f(g)=gvg^{-1}-v.$
\end{lemma}

 {\it Proof}.  If in \ref{theone} 
we set $g_1=x\in GL_n(F)$, $g_2=1+ \epsilon \xi\in GL_n(F[\epsilon])$ where $1$ is the identity,  $\epsilon^2=0$, and $\xi\in {\mathfrak gl}(F)$ and if we take the coefficient of $\epsilon$ 
we get
$$
d_x f(x\xi)=x d_1 f(\xi) x^{-1},
$$
equivalently
\begin{equation}\label{unu}
d_x f(\xi)=x d_1 f(x^{-1}\xi) x^{-1}.
\end{equation}
Similarly setting $g_1=1+\epsilon \xi$, $g_2=x$ in \ref{theone} we get
\begin{equation}\label{doi}
d_x f(\xi)= d_1 f(\xi x^{-1})+[\xi x^{-1},f(x)].
\end{equation}
From \ref{unu} and \ref{doi} we get
\begin{equation}\label{trei}
x d_1 f(x^{-1}\xi) x^{-1}=d_1 f(\xi x^{-1})+[\xi x^{-1},f(x)]
\end{equation}
for all $x\in GL_n(F)$ and $\xi\in {\mathfrak gl}_n(F)$ and hence for all $x\in GL_n(A)$, $\xi\in {\mathfrak gl}_n(A)$, $A$ an $F$-algebra. Setting in \ref{trei} $x=1+\epsilon \eta\in GL_n(F[\epsilon])$, with $\eta\in {\mathfrak gl}_n(F)$ and taking the coefficient of $\epsilon$ we get
$$d_1f([\xi,\eta])=[d_1f(\xi),\eta]+[\xi,d_1f(\eta)]$$
so $d_1f$ is a derivation on ${\mathfrak gl}_n(F)$. Note that $d_1f$ preserves ${\mathfrak sl}_n(F)$. Indeed by \ref{theone} we have
$$\text{tr}\  f(g_1g_2)=\text{tr} \ f(g_1)+\text{tr} \ f(g_2).$$
Since $GL_n$ has no regular non-zero homomorphism into ${\mathbb G}_a$ we get $\text{tr}\ f=0$ hence 
$f$ takes values in ${\mathfrak sl}_n(F)$. Now by the first Whitehead lemma \cite{jacobson}, p. 77, all derivations of ${\mathfrak sl}_n(F)$ are inner hence there exists $v\in {\mathfrak sl}_n(F)$ such that 
\begin{equation}
\label{patru}
d_1f(\xi)=[\xi,v]
\end{equation}
for all $\xi\in {\mathfrak sl}_n(F)$. By \ref{trei}, for all $\lambda \in F$, $f(\lambda 1)$ centralizes all of ${\mathfrak gl}_n(F)$ and, since $f(\lambda 1)\in{\mathfrak sl}_n(F)$, we must have $f(\lambda 1)=0$. We conclude that \ref{patru} holds for all $\xi\in {\mathfrak gl}_n(F)$. Let $f^v:GL_n(F)\ra {\mathfrak gl}_n(F)$ be the coboundary $f^v(g)=gvg^{-1}-v$ and note that $d_1f=d_1f^v$. Since, by \ref{trei}, $d_xf$ is determined by $d_1f$ it follows that $d_xf=d_xf^v$ for all $x\in GL_n(F)$. Since $f(1)=f^v(1)=0$ it follows that $f=f^v$.\qed

\medskip

As a consequence we get:

\begin{lemma}\label{cocoa2}
Assume we are in the  $\delta$-arithmetic setting.
Let $f:GL_n \ra {\mathfrak gl}_n$ be a  map of schemes over $R$  and assume $f$ is a cocycle  i.e. for all $g_1,g_2\in GL_n(R)$ equation \ref{theone} holds. Then 
 $f$ 
is a coboundary, i.e. there exists $v\in gl_n(R)$ such that for all $g\in GL_n(R)$,
 $f(g)=gvg^{-1}-v.$
\end{lemma}

{\it Proof}.
Let $K$ be the fraction field of $R$ and $F$ its algebraic closure. By Lemma \ref{cocoa} there exists $v_{F}\in {\mathfrak gl}_n(F)$ such that $f(g)=gv_{F}g^{-1}-v_{F}$ for all $g\in GL_n(F)$. Denoting by $v_K$  the sum of the conjugates of $v_{F}$ over $K$ 
we have
$f(g)=gv_Kg^{-1}-v_K$ for all $g\in GL_2(K)$. 
 Let $m\geq 1$ be an integer such that 
$p^{m}v_K\in {\mathfrak gl}_2(R)$.
  If  $w=p^{m}v_K-\alpha 1=(w_{ij})$ then we have
 $gwg^{-1}-w=p^{m}f(g)\in p^{m}{\mathfrak gl}_n(R)$
 for all $g\in GL_n(R)$. Taking $g$ to be an arbitrary diagonal matrix we get that $w_{ij}\equiv 0$ mod  $p^m$ for $i\neq j$. Taking $g=1+e_{ij}$, $i\neq j$, with $e_{ij}$ having all its entries $0$ except the $(i,j)$ entry which is $1$ we get $w_{ii}\equiv w_{jj}$ mod $p^m$; since $w_{11}\equiv 0$ mod $p^m$
 we get that all the entries of $w$ are divisible by $p^{m}$ in $R$.
 So $w=p^{m}v$, $v\in {\mathfrak gl}_n(R)$, and hence $f(g)=gvg^{-1}-v$ for all $g\in GL_n(R)$, which proves our claim.\qed

\begin{lemma}
\label{croco}
\

1) \cite{char} Assume we are in the $\delta$-arithmetic setting. Then any $\d$-homomorphism $f:{\mathbb G}_a\ra {\mathbb G}_a$, viewed as a map $f:R\ra R$,   has the form
 $$f(a)=\sum_{i=0}^r\lambda_i \phi^i(a),\ \ a\in R$$
 for some $\lambda_0,...,\lambda_r\in R$.
 
 2) \cite{cassidy} Assume we are in the  $\delta$-algebraic setting. Then any $\d$-homomorphism $f:{\mathbb G}_a\ra {\mathbb G}_a$, viewed as a map $f:R\ra R$,   has the form
 $$f(a)=\sum_{i=0}^r\lambda_i \d^i(a),\ \ a\in R$$
 for some $\lambda_0,...,\lambda_r\in R$.
 
\end{lemma}

{\it Proof.} 2) is due to Cassidy \cite{cassidy}, p. 936. To check 1), note that
by \cite{char} there exists $\nu\geq 1$ and  $\lambda_i\in R$ such that 
 $p^{\nu}f(a)=\sum\lambda_i \phi^i(a)$,  $a\in R$.
We are left  to check that if for $c_i\in R$  we have that for any $a\in R$, $\sum c_i \phi^i(a)\in pR$ then we must have $c_i\in pR$ for all $i$. But this follows from the fact that $\sum c_i a^{p^i}\equiv 0$ mod $p$ for all $a\in R$ and from the fact that the residue field of $R$ is algebraically closed.
\qed

\medskip

Next recall from \cite{char}, p. 313,  that  there exists a  remarkable $\d$-homomorphism $\psi:{\mathbb G}_m\ra {\mathbb G}_a$  given on $R$-points by
$$\psi(a)=\sum_{n\geq 1}(-1)^{n-1}\frac{p^{n-1}}{n}\left(\frac{\d a}{a^p}\right)^n.$$
With this notation we have:

\begin{lemma}
\label{allchar}
\

1) \cite{char} Assume we are in the  $\delta$-arithmetic setting. Then any $\d$-homomorphism $f:{\mathbb G}_m\ra {\mathbb G}_a$, viewed as a map $f:R^{\times}\ra R$,  has the form 
$$f(a)=\sum_{i=0}^r \lambda_i
\phi^i(\psi(a)),\ \ a\in R^{\times}$$ for some
 $\lambda_0,...,\lambda_r\in R$.
 
 2)  \cite{cassidy} Assume we are in the  $\delta$-algebraic setting. Then any $\d$-homomorphism $f:{\mathbb G}_m\ra {\mathbb G}_a$, viewed as a map $f:R^{\times}\ra R$,  has the form 
$$f(a)=\sum_{i=0}^r \lambda_i
\d^i(\d a \cdot a^{-1}),\ \ a\in R^{\times}$$ for some
 $\lambda_0,...,\lambda_r\in R$.

\end{lemma}

{\it Proof}. 2) is due to Cassidy \cite{cassidy}. To check 1) note that 
by \cite{char} there exist $\nu\geq 1$ such that $p^{\nu}f$ has the above form; we conclude exactly as in Lemma \ref{croco}.
\qed

\medskip

If instead of looking at $\d$-homomorphisms ${\mathbb G}_m\ra {\mathbb G}_a$ we are looking 
at cocycles of non-trivial actions the situation is very different:

\begin{lemma}
\label{aj}
Assume we are either in the $\delta$-algebraic setting or in the  $\delta$-arithmetic  setting.
Assume   $f:{\mathbb G}_m\ra {\mathbb G}_a$ is a $\d$-morphism and $0\neq s\in \bZ$ is an integer such that for any $a_1,a_2\in R^{\times}$ we have
$$f(a_1a_2)=f(a_1)+a_1^sf(a_2).$$
 Then there exists $\mu \in R$ such that 
$$f(a)=\mu(1-a^s),\ \ a\in R^{\times};$$
in particular $f$ is  a regular cocyle.
\end{lemma}

{\it Proof}.
Let us assume first we are in the $\delta$-arithmetic setting.
By induction on $r\geq 1$ one gets, for all $a\in R^{\times}$:
$$f(a^r)=(1+a^s+a^{2s}+...+a^{(r-1)s})f(a).$$
In particular for any $\nu\geq 0$ one gets
\begin{equation}
\label{cucu}
(1-a^s)f(a^{p^{\nu}})=(1-a^{sp^{\nu}})f(a)\end{equation}
Now recall that  
$$\cO(J^n({\mathbb G}_m))=R[x,x^{-1},x',...,x^{(n)}]\h.$$
Identify the point set map $f$ with a restricted power series
$f\in R[x,x^{-1},x',...,x^{(n)}]\h$ and write
$$f=\sum_{i_1,...,i_n\geq 0} f_{i_1...i_n} (x) (x')^{i_1}...(x^{(n)})^{i_n}$$
where $f_{i_1...i_n}\in R[x,x^{-1}]\h$ tend to $0$ $p$-adically for $|i_1|+...+|i_n|\ra \infty$.
By \ref{cucu} we get
\begin{equation}\label{iri}(1-x^s)\sum_{i_1,...,i_n\geq 0}
 f_{i_1...i_n} (x^{p^{\nu}}) \d(x^{p^{\nu}})^{i_1}...\d^n(x^{p^{\nu}})^{i_n}\end{equation}
$$=(1-x^{sp^{\nu}})\sum_{i_1,...,i_n\geq 0} f_{i_1...i_n} (x) (x')^{i_1}...(x^{(n)})^{i_n}$$
By induction on $\nu$ it is easy to see that  for any $\nu\geq r\geq 1$ we have 
\begin{equation}
\label{plan}
\d^r(x^{p^{\nu}})\in p^{\nu-r+1}R[x,x',...,x^{(r)}];\end{equation}
cf. also \cite{pfin1}.
Fix $(i_1,...,i_n)\neq (0,...,0)$ and let $\nu\geq n$ be arbitrary. Picking out the coefficient of $(x')^{i_1}...(x^{(n)})^{i_n}$ in the equation \ref{iri} we get
$$(1-x^{sp^{\nu}})f_{i_1...i_n}\in p^{\nu-n+1}R[x,x^{-1}]\h.$$
Since $\nu$ is arbitrary we get 
$$f_{i_1...i_n}=0 \ \ \text{for}\ \ \ (i_1,...,i_n)\neq (0,...,0),$$
hence $f=f_{0...0}$. For $\nu=1$ we may rewrite equation \ref{iri} as
$$(1-x^s)f_{0...0}(x^{p})=(1-x^{sp})f_{0...0}(x)$$
The latter is an equality in $R[x,x^{-1}]\h$ so we can also  view it as an equality in the ring
$$A=\{\sum_{i=-\infty}^{\infty} c_ix^i; c_i\ra 0\ \  \text{for}\ \ i\ra -\infty\}$$
If $s\geq 1$ we get the following equality in $A$:
$$g(x^{p})=g(x),\ \ \ g(x)=\frac{f_{0...0}(x)}{1-x^s}=f_{0...0}(x)(1+x^s+x^{2s}+...);$$
if $s\leq -1$ we get the following equality in $A$:
$$g(x^{p})=g(x),\ \ \ g(x)=\frac{-x^{-s}f_{0...0}(x)}{1-x^{-s}}=-x^{-s}f_{0...0}(x)(1+x^{-s}+x^{-2s}+...);$$
In either case we get $g=\mu\in R$ and we are done.

Assume now we are in the $\delta$-algebraic setting. Identify $f$ with a polynomial
in $A[x,x^{-1},x',...,x^{(n)}]$ where $A\subset R$ is a subring finitely generated over $\bZ$.
Now select a prime $p\in \bZ$ which is not invertible in $A$; hence by Krull's intersection theorem $\cap_{n=1}^{\infty} p^nA=0$. Then one can conclude by the same argument as above in which $R$ is replaced by $A$ (and all completions are dropped); by the way \ref{plan} is, in this case, obvious.
\qed

\begin{lemma}
\label{success}
Assume we are either in the $\delta$-algebraic setting or in the  $\delta$-arithmetic setting.
Assume  we are given a $\d$-map $\epsilon:{\mathbb G}_m\times {\mathbb G}_a^{n-1}\ra {\mathfrak gl}_{n-1}$ and  a row vector $\lambda\in R^{n-1}$ such that for any $(a_1,b_1), (a_2,b_2)\in R^{\times}\times R^{n-1}$:
 \begin{equation}
\label{train1}
\epsilon(a_1a_2,b_1+a_1b_2)  =  \epsilon(a_1,b_1)+\epsilon(a_2,b_2)+  a_1^{-1}(a_2^{-1}-1)\lambda^t b_1.
\end{equation}
Then, for any $a_1,a_2,a\in R^{\times}$, $b\in R^{n-1}$:
\begin{equation}
\label{camil1}
\epsilon(a_1a_2,0)=\epsilon(a_1,0)+\epsilon(a_2,0)
\end{equation}
\begin{equation}
\label{petrescu1}
\epsilon(a,b)=\epsilon(a,0)+ a^{-1}\lambda^t b.
\end{equation} \end{lemma}

{\it Proof}. Assume first we are in the  $\delta$-arithmetic setting.
Setting $b_1=b_2=0$ in equation \ref{train1} we get
\ref{camil1}. 
Setting $a_1=a_2=1$ in equation \ref{train1} we get
$$\epsilon(1,b_1+b_2)  = \epsilon(1,b_1)+\epsilon(1,b_2).
$$
So if $\epsilon=(\epsilon_{kl})$ then by Lemma \ref{croco} 
there exist $\lambda_{ikl}\in R^{n-1}$, $i\geq 0$, such that
\begin{equation}
\label{foc}
\epsilon_{kl}(1,b)=\sum \phi^i(b)\lambda_{ikl}^t,\ \ b \in R^{n-1}.\end{equation}
Now setting $a_1=a\in R^{\times}$, $a_2=1$, $b_1=0$, $b_2=a^{-1}b\in R^{n-1}$ in equation \ref{train1} we get
\begin{equation}
\label{trainy1}
\epsilon_{kl}(a,b)  = \epsilon_{kl}(a,0)+\epsilon_{kl}(1,a^{-1}b)=\epsilon_{kl}(a,0)+\sum  \phi^i(a)^{-1}\phi^i(b)\lambda_{ikl}^t.\end{equation}
Setting $a_1=1$, $a_2=a\in R^{\times}$, $b_1=b$, $b_2=0$, in equation \ref{train1} we get
\begin{equation}
\label{nice1}
\epsilon_{kl}(a,b)  = \epsilon_{kl}(1,b)+\epsilon_{kl}(a,0)+(a^{-1}-1)(\lambda^t b)_{kl}\end{equation}
hence, using equation \ref{trainy1} we get 
\begin{equation}
\label{mention1}
\sum (\phi^i(a)^{-1}-1)\phi^i(b)\lambda^t_{ikl}=(a^{-1}-1)(\lambda^tb)_{kl}.
\end{equation}
This can be viewed as an identity in $\phi^i(a),\phi^i(b)$ so we get $\lambda_{ilk}=0$ for all $i\geq 1$
and all $k,l$ and also $b\lambda_{0kl}^t=(\lambda^t b)_{kl}$. Then \ref{petrescu1} follows from \ref{trainy1}.

Assume now we are in the $\delta$-algebraic setting. The above argument can be modified as follows. Instead of equation \ref{foc} one gets
\begin{equation}
\label{focs}
\epsilon_{kl}(1,b)=\sum (\d^ib)\lambda_{ikl}^t,\ \ b \in R^{n-1}.\end{equation}
Instead of equation \ref{trainy1} one gets  
\begin{equation}
\label{trainy1s}
\epsilon_{kl}(a,b)  =\epsilon_{kl}(a,0)+\sum  \d^i(a^{-1}b)\lambda_{ikl}^t.\end{equation}
Equation \ref{nice1} remains then valid.
 Then instead of \ref{mention1} we get 
\begin{equation}
\label{mention1s}
\sum (\d^i(a^{-1}b)-\d^ib)\lambda^t_{ikl}=(a^{-1}-1)(\lambda^tb)_{kl}.
\end{equation}
Setting $b=e_j=(0,...,1,...,0)$ (with $1$ on the $j$th place) we get 
$$(a^{-1}-1)e_j\lambda^t_{0kl}+\sum_{i\geq 1} \d^i(a^{-1})e_j\lambda^t_{ikl}=(a^{-1}-1)\lambda_k \d_{lj}$$
($\d_{lj}$ the Kronecker symbol). Hence $e_j \lambda^t_{ikl}=0$ for $i\geq 1$ and all $j,k,l$.
Hence $\lambda_{ilk}=0$ for $i\geq 1$ and all $k,l$ and also $e_j\lambda^t_{0kl}=\lambda_k\d_{lj}$. So $a^{-1}b \lambda_{0kl}^t=(a^{-1}\lambda^tb)_{kl}$.
This ends the proof.
\qed
\bigskip

\begin{lemma}
\label{successs}
Assume we are either in the $\delta$-algebraic setting or in the $\delta$-arithmetic setting.
Assume  we are given a $\d$-map $\alpha:{\mathbb G}_m\times {\mathbb G}_a^{n-1}\ra {\mathbb A}^1$  and  a row vector $\lambda\in R^{n-1}$ such that for any $(a_1,b_1), (a_2,b_2)\in R^{\times}\times R^{n-1}$:
\begin{equation}
\label{asking}
\alpha(a_1a_2,b_1+a_1b_2)  =  \alpha(a_1,b_1)+\alpha(a_2,b_2)+  a_1^{-1}(1-a_2^{-1})b_1 \lambda^t.
\end{equation}
Then  for any $a_1,a_2,a\in R^{\times}$, $b\in R^{n-1}$:
\begin{equation}
\label{homo}
\alpha(a_1a_2,0)=\alpha(a_1,0)+\alpha(a_2,0)\end{equation}
\begin{equation}
\label{cap}
\alpha(a,b)=\alpha(a,0)- a^{-1}b\lambda^t.
\end{equation} 
\end{lemma}

{\it Proof}. Entirely similar to the proof of Lemma \ref{success}.\qed

\begin{lemma}
\label{successss}
Assume we are either in the $\delta$-algebraic setting or in the $\delta$-arithmetic setting.
Assume  we are given
a $\d$-map $\beta:{\mathbb G}_m\times {\mathbb G}_a^{n-1}\ra {\mathbb A}^{n-1}$, a $\d$-homomorphism $\psi:{\mathbb G}_m\ra {\mathfrak gl}_{n-1}$, and a row vector $\lambda\in R^{n-1}$ such that
for any $(a_1,b_1), (a_2,b_2)\in R^{\times}\times R^{n-1}$:
\begin{equation}
\label{maddd}
\begin{array}{rcl}
\beta(a_1a_2,b_1+a_1b_2) & =  &\beta(a_1,b_1)+a_1 \beta(a_2,b_2)\\
\  & \  & \  \\
\  & \  & +b_1 \psi(a_2)\\
\  & \  & \  \\
\  & \  & +
a_2^{-1}b_2\lambda^t b_1+a_2^{-1} b_1\lambda^tb_2+a_1^{-1}(a_2^{-1}-1)b_1 \lambda^t b_1\end{array}\end{equation}

1) Assume we are in the  $\delta$-arithmetic setting. Then $\psi=0$ and  there exist $\mu\in R^{n-1}$ and $\nu\in {\mathfrak gl}_{n-1}(R)$ such that for any $a\in R^{\times}$, $b\in R^{n-1}$ we have
\begin{equation}
\label{baby}
\beta(a,b)=(1-a)\mu+ a^{-1}b \lambda^t b+b\nu.
\end{equation}

2)  Assume we are in the  $\delta$-algebraic setting. Then   there exist $\mu\in R^{n-1}$ and  $\nu, \eta\in {\mathfrak gl}_{n-1}(R)$  such that for any $a\in R^{\times}$, $b\in R^{n-1}$ we have
\begin{equation}
\label{shoes}
\psi(a)=-(a^{-1}\d a)\nu;
\end{equation}
\begin{equation}
\label{baby}
\beta(a,b)=(1-a)\mu+ a^{-1}b \lambda^t b+b\eta+a\d(a^{-1}b)\nu.
\end{equation}

\end{lemma}

{\it Proof}. 
Assume first we are in the  $\delta$-arithmetic setting.
Setting $b_1=b_2=0$ in equation \ref{maddd} we get
$$\beta(a_1a_2,0)  = \beta(a_1,0)+a_1\beta(a_2,0).$$
So by Lemma \ref{aj} there exists $\mu\in R^{n-1}$ such that
$$\beta(a,0)=(1-a)\mu,\ \ a\in R^{\times}.$$
Setting $a_1=a_2=1$ in equation \ref{maddd} we get
$$\beta(1,b_1+b_2)  = \beta(1,b_1)+\beta(1,b_2)+b_2 \lambda^t b_1+b_1\lambda^t b_2,
$$
because $\psi(1)=0$. 
Let $\beta^*(b)=\beta(1,b)-b\lambda^t b$. Then
$$\beta^*(b_1+b_2)=\beta^*(b_1)+\beta^*(b_2).$$
So by Lemma \ref{croco}
there exist $\lambda_{i}\in {\mathfrak gl}_{n-1}(R)$, $i\geq 0$, such that
$$\beta^*(b)=\sum  \phi^i(b)\lambda_{i},\ \ b \in R^{n-1}$$
hence
\begin{equation}
\label{heat}
\beta(1,b)=b\lambda^t b+\sum  \phi^i(b)\lambda_i.\end{equation}
Now setting $a_1=a\in R^{\times}$, $a_2=1$, $b_1=0$, $b_2=a^{-1}b\in R^{n-1}$ in equation \ref{maddd} we get
\begin{equation}
\label{trainyyy}
\begin{array}{rcl}
\beta(a,b) & = & \beta(a,0)+a\beta(1,a^{-1}b)\\
\  & \  & \  \\
\  & = & (1-a)\mu+a\left( a^{-2}b \lambda^t b+\sum  \phi^i(a)^{-1}\phi^i(b)\lambda_i\right).\end{array}\end{equation}
Setting $a_1=1$, $a_2=a\in R^{\times}$, $b_1=b$, $b_2=0$ in equation \ref{maddd} we get
\begin{equation}
\label{nice}
\beta(a,b)  = \beta(1,b)+\beta(a,0)+b\psi(a)+(a^{-1}-1)b\lambda^tb,\end{equation}
hence, using equation \ref{trainyyy} we get, after cancellation: 
\begin{equation}
\label{mentionnn}
 \sum (a\phi^i(a)^{-1}-1)\phi^i(b)\lambda_i =b\psi(a).\end{equation}
Let $a=\zeta$ be an arbitrary root of unity. Since $\psi(\zeta)=0$ and $\phi(\zeta)=\zeta^p$ we get:
$$ \sum (\zeta^{1-p^i}-1)\phi^i(b)\lambda_i =0.$$
Since, for each $\zeta$ this can be seen as an identity in $\phi^i(b)$ we get that for all $\zeta$ and all $i\geq 0$,
$$ (\zeta^{1-p^i}-1)\lambda_i =0.$$
We get that $\lambda_i=0$ for all $i\geq 1$.
By \ref{mentionnn}  we get $\psi=0$. By \ref{trainyyy}, if $\nu:=\lambda_0$, we get \ref{baby}. This ends to the proof of assertion 1) in the Lemma.

Assume now we are in the $\delta$-algebraic setting. The above argument can be modified as follows. Instead of \ref{heat} we get
\begin{equation}
\label{heats}
\beta(1,b)=b\lambda^t b+\sum  (\d^i b)\lambda_i.\end{equation}
Instead of equation \ref{trainyyy} we get
\begin{equation}
\label{trainyyys}
\begin{array}{rcl}
\beta(a,b) = (1-a)\mu+a\left( a^{-2}b \lambda^t b+\sum  \d^i(a^{-1}b)\lambda_i\right).\end{array}\end{equation}
Equation \ref{nice} is valid.
 Instead of equation \ref{mentionnn} we then get
\begin{equation}
\label{mentionnns}
 \sum (a\d^i(a^{-1}b)-\d^i b)\lambda_i =b\psi(a).\end{equation}
Let $\psi=(\psi_{kl})$ and $\lambda_i=(\lambda_{ikl})$. Setting $b=e_j=(0,...,1,...,0)$ we get
$$\sum_{i\geq 1} a\d^i(a^{-1})\lambda_{ijl}=\psi_{jl}(a)$$
Replacing $a$ by $a^{-1}$ and using $\psi(a^{-1})=-\psi(a)$ we get
$$\sum_{i\geq 1} (\d^ia)\lambda_{ijl}=-a\psi_{jl}(a)$$
Using the fact that  $\psi_{jl}$ is a $\d$-homomorphism ${\mathbb G}_m\ra {\mathbb G}_a$ and $\d$ is a $\d$-homomorphism ${\mathbb G}_a\ra {\mathbb G}_a$  and setting $a=a_1a_2$ one immediately concludes that
$\lambda_{ijl}=0$ for $i\geq 1$ and all $j,l$. Hence $\psi(a)=-(a^{-1}\d a)\lambda_1$.
Setting $\nu=\lambda_1$ and $\eta=\lambda_0$ we conclude by \ref{trainyyys}.
 \qed

\begin{lemma}
\label{indet}
Let $x$ be an $n\times n$ matrix of indeterminates over $\bZ$ and let 
$\Delta_i(x)$ be the determinant of the matrix obtained from $x$ by removing
the first $i$ rows and the first $i$ columns. Let $\Delta(x)=\prod_{i=1}^{n-1}\Delta_i(x)$.  Let
$W$ be the group of all $n\times n$ matrices obtained from the identity matrix $1_n$ by permuting its columns. Then there exist
$w_0,...,w_N \in W$ and there exist elements
$$a_1(x),...,a_N(x)\in \bZ[x,\Delta(x)^{-1}],$$
and $1\times (n-1)$ vectors
$$b_1(x),...,b_N(x)\in \bZ[x,\Delta(x)^{-1}]^{n-1},$$
such that if
$$s_i(x)=\left(\begin{array}{ll} a_i(x) & b_i(x)\\ 0 & 1_{n-1}\end{array}\right)$$
then
$$x=w_0 \cdot s_1(x)\cdot  w_1 \cdot s_2(x)\cdot  w_2 \cdot ... \cdot s_N(x)\cdot w_N$$
\end{lemma}

{\it Proof}.
We proceed by induction on $n$. Write 
$$x=\left(\begin{array}{ll} u & y\\ z^t & w\end{array}\right)$$
where $u$ is one variable, $y,z$ are $1\times (n-1)$ matrices of variables and $w$ is an $(n-1)\times (n-1)$ matrix of variables. We may write
\begin{equation}
\label{foamee}
x=\left(\begin{array}{ll} u-yw^{-1}z^t & yw^{-1}\\ 0& 1_{n-1}\end{array}\right)
\left(\begin{array}{ll} 1 & 0\\ 0 & w\end{array}\right) 
\left(\begin{array}{ll} 1 & 0\\ w^{-1}z^t & 1_{n-1}\end{array}\right).
\end{equation}
Now for $n \times n$ matrices $m,m'$ we say $m$ and $m'$ are $W$-equivalent if $W mW=W m' W$.
We then have that $\left(\begin{array}{ll} 1 & 0 \\ 0 & w\end{array}\right)$ is $W$-equivalent to
$\left(\begin{array}{ll} w & 0\\ 0 & 1\end{array}\right)$.
Also, if $w^{-1}z^t=(c_1(x),...,c_{n-1}(x))^t$ then
$$\left(\begin{array}{ll} 1 & 0\\ w^{-1}z^t & 1_{n-1}\end{array}\right)=\left(\begin{array}{lll} 1 & 0 & 0 \\ c_1(x) & 1 & 0\\
0 & 0 & 1_{n-2}\end{array}\right)\left(\begin{array}{llll} 1 & 0 & 0 & 0\\  0 & 1 & 0 & 0\\
c_2(x) & 0 & 1 & 0\\
0 & 0 & 0 &  1_{n-3}\end{array}\right)....
$$
and each of the matrices in the right hand side of the above equation is $W$-equivalent to its transpose.
We may conclude from \ref{foamee} by the induction hypothesis applied to $w$.
\qed

\bigskip

{\it Proof of Theorem \ref{main}}. 
Write matrices in $GL_{n}(R)$ in block from $\left(\begin{array}{ll}a & b\\ c^t & d\end{array}\right)$ where $a\in R$, $b,c\in R^{n-1}$ are viewed as a $1\times (n-1)$ matrices, $c^t$ is the transpose of $c$, and $d$ is an $(n-1)\times (n-1)$ matrix. Also let $1_{n-1}$ be the identity $(n-1)\times (n-1)$ matrix.
 We may write
$$f\left(\begin{array}{ll}a & b\\ 0 & 1_{n-1}\end{array}\right)=
\left(\begin{array}{cc}\alpha(a,b) & \beta(a,b)\\ \gamma(a,b)^t & \epsilon(a,b)\end{array}\right)$$
where $\alpha, \beta,\gamma, \epsilon$ are $\d$-morphisms from ${\mathbb G}_m\times {\mathbb G}_a^{n-1}$ to ${\mathbb A}^1,{\mathbb A}^{n-1},{\mathbb A}^{n-1}, {\mathfrak gl}_{n-1}$ respectively.
If we let $(a_i,b_i)\in R^{\times}\times R^{n-1}$, $i=1,2$, and if we set 
$$
\alpha_i=\alpha(a_i,b_i), \ \ \beta_i=\beta(a_i,b_i), \ \ \gamma_i=\gamma(a_i,b_i),\ \  \epsilon_i=\epsilon(a_i,b_i),\ \ i=1,2,
$$
then the cocycle relation for $f$ yields:
\begin{equation}
\label{list}
\begin{array}{rrcl}
(1) & \alpha(a_1a_2,b_1+a_1b_2) & = & \alpha_1+\alpha_2+a_1^{-1}b_1\gamma_2^t\\
(2) & \beta(a_1a_2,b_1+a_1b_2) & = & \beta_1+a_1 \beta_2-\alpha_2 b_1+b_1 \epsilon_2 -a_1^{-1}b_1\gamma_2^t b_1\\
(3) & \gamma(a_1a_2,b_1+a_1b_2) & = & \gamma_1+a_1^{-1}\gamma_2\\
(4) & \epsilon(a_1a_2,b_1+a_1b_2) & = & \epsilon_1+\epsilon_2-a_1^{-1}\gamma_2^t b_1
\end{array}
\end{equation}
By equation \ref{list} (3) and by Lemma \ref{aj} there exists $\lambda\in R^{n-1}$ such that
\begin{equation}
\label{official}
\gamma(a,b)=(1-a^{-1})\lambda,\ \ a\in R^{\times}, b\in R^{n-1}.
\end{equation}
So equation \ref{list} (4) reads:
\begin{equation}
\label{ask}
\epsilon(a_1a_2,b_1+a_1b_2)  =  \epsilon(a_1,b_1)+\epsilon(a_2,b_2)+  a_1^{-1}(a_2^{-1}-1)\lambda^t b_1.
\end{equation}
By Lemma 
\ref{success}
we get 
\begin{equation}
\label{homoe}
\epsilon(a_1a_2,0)=\epsilon(a_1,0)+\epsilon(a_2,0)\end{equation}
\begin{equation}
\label{cape}
\epsilon(a,b)=\epsilon(a,0)+ a^{-1}\lambda^t b.
\end{equation} 
Similarly,  equation \ref{list} (1) reads
\begin{equation}
\label{asking}
\alpha(a_1a_2,b_1+a_1b_2)  =  \alpha(a_1,b_1)+\alpha(a_2,b_2)+  a_1^{-1}(1-a_2^{-1})b_1 \lambda^t.
\end{equation}
By Lemma 
\ref{successs}
we get
\begin{equation}
\label{homo}
\alpha(a_1a_2,0)=\alpha(a_1,0)+\alpha(a_2,0)\end{equation}
\begin{equation}
\label{cap}
\alpha(a,b)=\alpha(a,0)- a^{-1}b\lambda^t.
\end{equation} 
Finally by equations \ref{official}, \ref{cape},  \ref{cap},  equation \ref{list} (2) reads
\begin{equation}
\label{mad}
\begin{array}{rcl}
\beta(a_1a_2,b_1+a_1b_2)  & = & \beta(a_1,b_1)+a_1 \beta(a_2,b_2) \\
\  & \  & \  \\
 \ & \  & +
b_1(\epsilon(a_2,0)-\alpha(a_2,0)1_{n-1})\\
\  & \  & \  \\
\  & \  & +
a_2^{-1}b_2\lambda^t b_1+a_2^{-1} b_1\lambda^tb_2+a_1^{-1}(a_2^{-1}-1)b_1 \lambda^t b_1.\end{array}\end{equation}
By Lemma \ref{successss} 
\begin{equation}
\label{zgomot}
\epsilon(a,0)=\alpha(a,0)1_{n-1}, \ a\in R^{\times}\end{equation}
and 
there exist $\mu\in R^{n-1}$, $\nu\in {\mathfrak gl}_{n-1}(R)$ such that
\begin{equation}
\label{shack}
\beta(a,b)=(1-a)\mu+ a^{-1}b \lambda^t b+b\nu.
\end{equation}
Let $\omega:{\mathbb G}_m\ra {\mathbb G}_a$ be the $\d$-homomorphism defined by $\omega(a)=\alpha(a,0)$.
Consider the classical $\d$-cocycle $f_0:GL_n\ra {\mathfrak gl}_n$ defined by
$$f_0(g):=f(g)-\omega(\det\ g)1_n.$$
By Lemma \ref{cocoa2} it is sufficient to show that $f_0$ is induced on points by a regular map. 
Let $H\subset GL_n$ the natural closed subgroup scheme such that
$$H(R)=\{\left(\begin{array}{ll}a & b \\ 0 & 1_{n-1}\end{array}\right);a\in R^{\times}, b\in R^{n-1}\}.$$
By equations  \ref{cap}, \ref{cape}, \ref{official}, the restriction of $f_0$ to $H$  is a regular map. On the other hand by Lemma \ref{indet} there exist finitely  many matrices $w_0,w_1,...,w_N\in W$ such that  the map
$$H\times...\times H\ra GL_n$$ defined on points by
$$ (h_1,...,h_N)\mapsto w_0 h_1w_1 \cdot\cdot\cdot h_N w_N$$
has a rational section 
$$s:U\ra H\times ... \times H, \ s(g)=(s_1(g),...,s_N(g))$$
defined on
 a Zariski open set $U$ of $GL_n$ which meets the special fiber at $p$. So for 
 $g\in U(R)$ we have
 $$f_0(g)=f_0(w_0\cdot  s_1(g)\cdot w _1 \cdot...\cdot s_N(g)\cdot  w_N).$$
 By the cocycle condition we get that the restriction of $f_0$ to $U$ is regular and hence, using translations and the cocycle condition again we get that $f_0$ is regular on the whole of $GL_n$. \qed
   
 \bigskip
 
 {\it Proof of Theorem \ref{secondary}}.
 One can redo the argument in the proof of Theorem \ref{main} as follows. 
 First we assume $f$ is any classical $\d$-cocycle (not necessarily $\d$-coherent).
 The formulae giving $\gamma, \alpha, \epsilon$
 (\ref{official}, \ref{homoe}, \ref{cape}, \ref{homo}, \ref{cap})
  are unchanged. Hence equation \ref{mad} is valid.
 By Lemma \ref{successss} there exist $\mu,\lambda \in R^{n-1}$ and  $\nu,\eta\in {\mathfrak gl}_{n-1}$ such that
\begin{equation}
\label{zgomot}
\epsilon(a,0)=\alpha(a,0)1_{n-1}-(a^{-1}\d a)\nu, \ a\in R^{\times}\end{equation}
\begin{equation}
\label{shack}
\beta(a,b)=(1-a)\mu+ a^{-1}b \lambda^t b+b\eta+a\d(a^{-1}b)\nu.
\end{equation}
Set  $g=\left(\begin{array}{ll} a & b \\ 0 & 1_{n-1}\end{array}\right)\in H(R)$. By the above we conclude that 
\begin{equation}
\label{trip}
f(g)= \left(\begin{array}{lll}
\alpha(a,0)1_n-a^{-1}b\lambda^t & (1-a)\mu+ a^{-1}b \lambda^t b+b\eta+a\d(a^{-1}b)\nu\\
\  &\  & \  \\
(1-a^{-1})\lambda^t & \alpha(a,0)1_{n-1}-(a^{-1}\d a)\nu+a^{-1}\lambda^t b\end{array}\right)
\end{equation}
Now let us assume $f$ is $\d$-coherent. Since $f$ must send $H$ into its Lie algebra, consisting of the matrices $\left(\begin{array}{ll} \star & \star \\ 0 & 0_{n-1}\end{array}\right)$, we get that $\lambda=0$ and $$\alpha(a,0)1_{n-1}=(a^{-1}\d a)\nu.$$
This forces $\nu$ to be a scalar matrix, say $\nu=\epsilon 1_{n-1}$ and $\alpha(a,0)=\epsilon \cdot a^{-1}\d a$. 
 If for $v\in {\mathfrak gl}_n(R)$ we write $f^v(g)=gvg^{-1}-v$, $g\in GL_n(R)$, then note that
 $$f^{v}\left(\begin{array}{ll} a & b \\0 & 1_{n-1}\end{array}\right)=\left(\begin{array}{cc} 0 & (1-a)\mu +b\eta \\ 0 & 0 \end{array}\right)\ \ \text{for}\ \ v=\left(\begin{array}{cc} 0 & -\mu\\0 & \eta\end{array}\right)$$
We conclude that 
$$f(g)=f^v(g)+\epsilon \cdot \d g \cdot g^{-1}$$
for all $g\in H(R)$ and hence (as in the proof of Theorem \ref{main}) for all $g\in GL_n(R)$.
Now since $f$ and $g\mapsto \d g \cdot g^{-1}$ are $\d$-coherent it follows that $f^v$ is $\d$-coherent. It is enough to prove that $v$ is a scalar matrix. Assume $v$ is not scalar. Then we claim there exists $u\in GL_n(R^{\d})$ such that $uvu^{-1}$ is not in the Lie algebra  ${\mathfrak t}(R)\subset {\mathfrak gl}_n(R)$ of diagonal matrices. Indeed if one assumes
$$\{uvu^{-1};u\in GL_n(R^{\d})\}\subset {\mathfrak t}(R),$$
since $\{uvu^{-1};u\in GL_n(R^{\d})\}$ is Zariski dense in $\{uvu^{-1};u\in GL_n(R)\}$
we get that $\{uvu^{-1};u\in GL_n(R)\}\subset {\mathfrak t}(R)$. Note that $v$ itself is in ${\mathfrak t}(R)$. So the stabilizer of $v$ in $GL_n(R)$ under the adjoint action has dimension
$n_1^2+...+n_s^2$ where $n_1,...,n_s$ are the multiplicities of the eigenvalues of $v$;
note that $s\geq 2$ and $n_1+...+n_s=n$. Hence the dimension of $\{uvu^{-1};u\in GL_n(R)\}$ is $n^2-(n_1^2+...+n_s^2)> n=\dim {\mathfrak t}$ (because $s\geq 2$), a contradiction. 
Let now $G=u^{-1}Tu$ where $T\subset GL_n$ is the torus of diagonal matrices; so $G$ is defined over $R^{\d}$ and has Lie algebra $L(G)=u^{-1}L(T)u$.
Since $f^v$ is $\d$-regular we have $u^{-1}duvu^{-1}d^{-1}u-v\in u^{-1}L(T)(R)u$ for all $d\in T(R)$. Hence $duvu^{-1}d^{-1} -uvu^{-1}\in {\mathfrak t}$ for all $d\in T(R)$.
But $uvu^{-1}$ is not diagonal so there are indices $i\neq j$ such that then $(i,j)$-entry of $uvu^{-1}$ is non zero. Taking $d$ with diagonal entries $d_i\neq d_j$ we get
that the $(i,j)$-entry of $duvu^{-1}d^{-1} -uvu^{-1}$ is non-zero, a contradiction. This ends the proof.
\qed


\begin{thebibliography}{AA}

\bibitem{borger} J. Borger, {\it The basic geometry of Witt vectors, I}:
the affine case,  Algebra and Number Theory 5 (2011), no. 2, pp 231-285. 

\bibitem{borger2}
 Borger, J., {\it The basic geometry of Witt vectors, II: Spaces},
arXiv:1006.0092v1.

\bibitem{borgerf1} J. Borger, {\it $\Lambda$-rings and the field with one element},
 	arXiv:0906.3146 [math.NT]
  	
\bibitem{annals} A. Buium, {\it Intersections in jet spaces and a
conjecture of S.Lang}, Ann. of Math. 136 (1992), pp. 583-593.

\bibitem{lnm} A. Buium, {\it Differential algebraic groups of finite dimension}, LNM 1506, Springer, 1992.





\bibitem{hermann} A. Buium, {\it Differential Algebra and Diophantine
Geometry}, Hermann, Paris, 1994.


\bibitem{char}
A. Buium, {\it Differential characters of Abelian varieties over
$p$-adic fields}, Invent. Math., 122 (1995), pp. 309-340.


\bibitem {pjets}
A. Buium, {\it Geometry of $p$-jets}, Duke J. Math. 82, (1996), 2, pp. 349-367.

\bibitem{simple} A. Buium,  {\it Differential subgroups of simple algebraic groups over p-adic
fields},
Amer. J. Math. 120 (1998), 1277-1287.


\bibitem{book}
A. Buium, {\it Arithmetic Differential Equations}, Math. Surveys and
Monographs, 118, American Mathematical Society, Providence, RI, 2005.
xxxii+310 pp.


\bibitem{pfin1} A. Buium, {\it $p$-jets of finite algebras, I: $p$-divisible groups}, Documenta Mathematica, to appear.

\bibitem{adel2} A. Buium, T. Dupuy, {\it Arithmetic differential equations on $GL_n$, II: arithmetic Lie theory}, preprint.

\bibitem{adel3} A. Buium, T. Dupuy, {\it Arithmetic differential equations on $GL_n$, III: Galois groups}, preprint.

\bibitem{cassidy} P. J. Cassidy, {\it Differential algebraic groups}, Amer. J. Math. 94 (1972), 891-954.

\bibitem{cassidy2} P. J. Cassidy, {\it The classification of semisimple differential algebraic groups and the linear semisimple differential algebraic Lie algebras}, J. Algebra (1989), 169-238.

\bibitem{jacobson} N. Jacobson, {\it Lie Algebras}, Dover, New York, 1979.

\bibitem{kolchin}
E.R. Kolchin, {\it Differential algebra and algebraic groups}. Pure and
Applied Mathematics, Vol. 54. Academic Press, New York-London, 1973.
xviii+446 pp.

\bibitem{manin}
Yu.I. Manin, {\it Algebraic curves over fields with differentiation}, Izv.
Akad. Nauk SSSR, Ser. Mat. 22 (1958), pp. 737-756 = AMS translations Series 2,
37 (1964), pp. 59-78.

\bibitem{pillay} A. Pillay, {\it Algebraic $D$-groups and differential Galois theory}, Pacific J. Math 216,
2004, 343-360.


\end{thebibliography}
\end{document}